\documentclass[production,11pt]{jsg}\usepackage{latexsym}\usepackage{amsmath}\usepackage{xspace}\usepackage{enumerate}\usepackage{amsfonts}\usepackage{amscd}\usepackage{syntonly}\usepackage{amssymb}

\newtheorem{thm}{Theorem}
\newtheorem{prop}[thm]{Proposition}
\newtheorem{lemma}[thm]{Lemma}
\newtheorem{claim}{Claim}
\newtheorem{defi}[thm]{Definition}
\newtheorem{rmk}[thm]{Remark}
\newtheorem{cor}[thm]{Corollary}

\newenvironment{enui}{\begin{enumerate}[(i)]}{\end{enumerate}}

\def\nn{{\nonumber}}

\newcommand\wt[1]{{\widetilde{#1}}}
\newcommand{\BAR}[1]{{\overline{#1}}}

\def\ga{\gamma}
\def\de{\delta}
\def\eps{\varepsilon}

\def\om{\omega}
\def\lam{\lambda}
\def\Si{\Sigma}

\def\ze{\zeta}
\renewcommand\phi{\varphi}

\newcommand{\N}{\mathbb{N}}
\newcommand{\Z}{\mathbb{Z}}

\newcommand{\R}{\mathbb{R}}
\def\C{\mathbb C}

\def\D{\mathbb{D}}

\def\CP{\C{\operatorname{P}}}
\def\G{{\operatorname{G}}}

\def\g{\mathfrak g}
\def\A{\mathcal A}
\def\Lie{\operatorname{Lie}}
\def\Crit{\operatorname{Crit}}
\def\Pr{\operatorname{Pr}}

\def\sub{\subseteq}

\def\x{\times}
\def\wo{\setminus}

\def\Then{\quad\Longrightarrow\quad}

\def\id{{\operatorname{id}}}
\def\one{\mathbf{1}}

\def\iso{\cong}
\def\const{\equiv}
\def\d{\partial}
\def\na{\nabla}
\def\La{\Delta}
\def\lan{\langle}
\def\ran{\rangle}

\def\cf{cf. }
\def\ie{i.e.\xspace}

\def\wrt{w.r.t.\xspace}
\def\etc{etc.\xspace}
\def\Wlog{w.l.o.g.\xspace}

\def\loc{{\operatorname{loc}}}

\def\op{{\operatorname{op}}}

\title{The Invariant Symplectic Action and Decay for Vortices} 
\author{Fabian Ziltener (University of Toronto)}
\begin{document}
\begin{abstract} The (local) invariant symplectic action functional $\A$ is associated to a Hamiltonian action of a compact connected Lie group $\G$ on a symplectic manifold $(M,\om)$, endowed with a $\G$-invariant Riemannian metric $\lan\cdot,\cdot\ran_M$. It is defined on the set of pairs of loops $(x,\xi):S^1\to M\x\Lie\G$ for which $x$ satisfies some admissibility condition. I prove a sharp isoperimetric inequality for $\A$ if $\lan\cdot,\cdot\ran_M$ is induced by some $\om$-compatible and $\G$-invariant almost complex structure $J$, and, as an application, an optimal result about the decay at $\infty$ of symplectic vortices on the half-cylinder $[0,\infty)\x S^1$.
\end{abstract}
\maketitle
\tableofcontents
\section{Motivation and main results}\label{sec:main}
Let $(M,\om)$ be a symplectic manifold without boundary, and $\G$ be a compact connected Lie group with Lie algebra $\g$. Suppose that $\G$ acts on $M$ in a Hamiltonian way, with (equivariant) moment map $\mu:M\to\g^*$. We denote by $\lan\cdot,\cdot\ran:\g^*\x\g\to\R$ the natural contraction. Furthermore, we fix a $\G$-invariant metric $\lan\cdot,\cdot\ran_M$ on $M$, and denote by $|v|$, $\iota_x$ and $\iota_X:=\inf_{x\in X}\iota_x\geq0$ the corresponding norm of a vector $v\in TM$ and the injectivity radius of a point $x\in M$ and a subset $X\sub M$ respectively. For a smooth loop $x:S^1\to M$ of length $\ell(x)$ less than $2\iota_{x(S^1)}$ we denote by $\A(x)$ its (usual) symplectic action (see Section \ref{subsec:back}). We identify $S^1\iso\R/\Z$ and call a loop $x\in C^\infty(S^1,M)$ \emph{admissible} iff there exists a gauge transformation $g\in C^\infty(S^1,\G)$ such that $\ell(gx)<2\iota_{x(S^1)},$ and
\begin{equation}\label{eq:A wt g x int}\A(\wt gx)-\A(gx)=\int_0^1\big\lan\mu\circ x,\wt g^{-1}\dot{\wt g}-g^{-1}\dot g\big\ran\,dt,
\end{equation}
for every $\wt g\in C^\infty(S^1,\G)$ satisfying $\ell(\wt gx)\leq \ell(gx)$. 
\begin{defi}\label{defi:action} Let $(x,\xi)\in C^\infty(S^1,M\x\g)$ be a pair of loops, such that $x$ is admissible. We define the \emph{invariant (symplectic) action of $(x,\xi)$} to be 
\begin{equation}\label{eq:A x xi}\A(x,\xi):=\A(gx)+\int_0^1\big\lan\mu\circ x,\xi-g^{-1}\dot g\big\ran\,dt,\end{equation}
where $g\in C^\infty(S^1,\G)$ is chosen as above.
\end{defi}
This is a modified version of the ``local equivariant symplectic action functional'' introduced by A. R. Gaio and D. A. Salamon in \cite{GS}. More precisely, for $x\in M$ we denote by $L_x:\g\to T_xM$ the infinitesimal action of the Lie algebra on the tangent space to $M$ at $x$. Furthermore, we fix a $\G$-invariant inner product $\lan\cdot,\cdot\ran_\g$ on $\g$ and denote $|\xi|:=\sqrt{\lan\xi,\xi\ran_\g}$, for $\xi\in\g$. The induced operator norm on $\g^*$ is denoted by $|\cdot|_\op$. Gaio and Salamon define the action for pairs $(x,\xi)$ for which $|\mu\circ x|_\op$ and the twisted length $\ell(x,\xi):=\int_0^1|\dot x+L_x\xi|\,dt$ are small.\footnote{There is a gap in that definition, since the imposed smallness conditions do actually not guarantee that the action of $(x,\xi)$ is well-defined, see section \ref{subsec:back}. In that subsection a more direct way of fixing the gap is also mentioned.} To formulate the first main result of this paper, we denote by $M^*\sub M$ the subset of all points on which $\G$ acts freely, and by $\G x\in M^*/\G$ the orbit of a point $x\in M^*$. For a loop $\bar x:S^1\to M^*/\G$ we denote by $\bar\ell(\bar x)$ its length \wrt the Riemannian metric on $M^*/\G$ induced by $\lan\cdot,\cdot\ran_M$. Furthermore, for each subset $X\sub M$ we define 
\begin{equation}\label{eq:m X}m_X:=\inf\big\{|L_x\xi|\,\big|\,x\in X,\,\xi\in\g:\,|\xi|=1\big\}.
\end{equation}
For $p\in[1,\infty]$ and smooth loops $v:S^1\to TM$ and $\phi:S^1\to \g^*$ we denote by $||v||_p$ and $||\phi||_p$ the $L^p$-norms \wrt the Haar measure on $S^1$, the metric $\lan\cdot,\cdot\ran_M$ and the norm $|\cdot|_\op$ on $\g^*$. 
\begin{thm}[Sharp isoperimetric inequality]\label{thm:isoperi} Assume that there exists a $\G$-invariant $\om$-compatible almost complex structure $J$ such that $\lan\cdot,\cdot\ran_M=g_{\om,J}=\om(\cdot,J\cdot)$. Then for every compact subset $K\sub M^*$ and every constant $c>\frac1{4\pi}$ there exists a constant $\de>0$ with the following property. Suppose that $x\in C^\infty(S^1,K)$ is a loop satisfying $\bar\ell(\G x)<\de.$ Then $x$ is admissible, and for every loop $\xi\in C^\infty(S^1,\g)$ and every number $1\leq p\leq2$ we have
\begin{equation}
  \label{eq:isoperi}|\A(x,\xi)|\leq c||\dot x+L_x\xi||_p^2+\frac\pi{m_K^2}||\mu\circ x||_{\frac p{p-1}}^2.
\end{equation}
\end{thm}
This result generalizes the isoperimetric inequality for the usual action (\cf Theorem 4.4.1 in the book \cite{MS} by D. McDuff and D. A. Salamon), which corresponds to the case of the trivial Lie group $\G:=\{\one\}$. It is sharp in the sense that in general, the statement for $c<1/(4\pi)$ is wrong. An example illustrating this is given by $\G:=\{\one\}$ and the plane $M:=\R^2$ with the standard symplectic and Riemannian structures $\om_0$ and $\lan\cdot,\cdot\ran_0$. The constant $\pi/m_K^2$ is also optimal. To see this, consider 
\[M:=\R^2,\quad \om:=\om_0,\quad \lan\cdot,\cdot\ran_M:=\lan\cdot,\cdot\ran_0,\quad K:=\{(2,0)\},\] 
and the action of $\G:=S^1\iso\R/\Z$ on $\R^2$ by rotation, with moment map 
\begin{equation}\label{eq:mu R 2}\mu:\R^2\to \Lie(S^1)^*\iso\R^*,\quad \lan\mu(x),\xi\ran:=\pi \xi\big(1-|x|^2\big).\end{equation}
Let $C<\pi/m_K^2=1/(16\pi)$. Then the constant pair $(x,\xi):=\big((2,0),3/8\big)$ violates inequality (\ref{eq:isoperi}) with $\pi/m_K^2$ replaced by $C$, for every $p$, if we choose $c\in\big(1/(4\pi),1/(2\pi)-4C\big)$. 

The proof of Theorem \ref{thm:isoperi} is based on the isoperimetric inequality for the usual symplectic action and on an estimate for the holonomy of a connection around a loop in terms of the curvature of the connection and the length of the loop. Note that Gaio and Salamon proved an isoperimetric inequality for their equivariant action, for $p=2$ and a \emph{large} constant, \cf Lemma 11.3 in \cite{GS}.

To explain the application of Theorem \ref{thm:isoperi}, let $J$ be a $\G$-invariant and $\om$-compatible almost complex structure on $M$, $(\Si,j)$ be a Riemann surface equipped with a compatible area form $\om_\Si$, and let $P$ be a principal $\G$-bundle over $\Si$. The (symplectic) vortex equations for a pair $(u,A)$ are given by
\begin{equation}
  \label{eq:vort P}\left\{
  \begin{array}{ccc}
\bar \d_{J,A}(u)&=&0\\
F_A+(\mu\circ u)\om_\Si&=&0.
  \end{array}\right.
\end{equation}
Here $u$ is an equivariant map from $P$ to $M$, and $A$ is a connection one-form on $P$. Furthermore, $\bar\d_{J,A}(u)$ denotes the complex anti-linear part of $d_Au:=du+L_uA$, which we think of as a one-form on $\Si$ with values in the complex vector bundle $u^*TM/\G\to\Si$. Similarly, we view the curvature $F_A$ of $A$ as a two-form on $\Si$ with values in the adjoint bundle $\g_P:=(P\x\g)/\G\to\Si$. Finally, we identify $\g^*$ with $\g$ via $\lan\cdot,\cdot\ran_\g$, and we view $\mu\circ u$ as a section of $\g_P$. The equations (\ref{eq:vort P}) were discovered, independently, on the one hand by K. Cieliebak, A. R. Gaio and D. A. Salamon \cite{CGS}, and on the other hand by I. Mundet i Riera \cite{Mu1}, \cite{Mu2}. 

We fix now a number $a>0$ and consider the case in which the Riemann surface is the half-cylinder $\Si:=\big\{s+it\in \C\,\big|\,s\geq0\big\}/ai\Z,$ endowed with the standard complex structure $j:=i$. We denote $\La:=\d_s^2+\d_t^2$ and define $m_{\mu^{-1}(0)}$ as in (\ref{eq:m X}). We call an area form $\om_\Si=\lam^2ds\wedge dt$ on $\Si$ \emph{admissible} iff
\begin{equation}\label{eq:lam sup}\lam\geq \displaystyle\frac{2\pi}{am_{\mu^{-1}(0)}},\qquad \displaystyle\sup_{\Si}\big(|d(\lam^{-1})|^2+\La(\lam^{-2})\big)<2m_{\mu^{-1}(0)}^2.\end{equation}
Consider the following hypothesis.\\

\noindent\label{H}{\bf(H)} There exists a number $\eps>0$ such that the set $\big\{x\in M\,\big|\,|\mu(x)|\leq\eps\big\}$ is compact. Furthermore, $\G$ acts freely on $\mu^{-1}(0)$. \\

We fix a pair $w:=(u,A)$ where $u$ is an equivariant map from $P$ to $M$ and $A$ is a connection one-form on $P$. Recall that the energy density and energy of $w$ are defined by 
\begin{equation}\label{eq:e w E}e_w:=\frac12\Big(|d_Au|^2+|F_A|^2+|\mu\circ u|^2\Big),\qquad E(w):=\int_\Si e_w\om_\Si,\end{equation}
with the norms taken \wrt the metrics $g_{\om_\Si,j}:=\om_\Si(\cdot,j\cdot)$ on $\Si$ and $g_{\om,J}$ on $M$. The application of Theorem \ref{thm:isoperi} is the following.
\begin{thm}[Optimal decay for vortices on the half-cylinder]\label{thm:decay} Assume that hypothesis {\rm (H)} is satisfied. Let $\Si$ be the half-cylinder, $\om_\Si$ be an admissible area form on $\Si$, and $p>2$. Assume that $w:=(u,A)$ is a locally $W^{1,p}$-solution of the equations {\rm (\ref{eq:vort P})}, such that $E(w)<\infty$ and the image of $u$ has compact closure in $M$. Then for every $\eps>0$ there exists a constant $C$ such that
\begin{equation} \label{eq:d A u lam} e_w(s+it)\leq C\lam^{-2}e^{\left(-\frac{4\pi}a+\eps\right)s},\quad \forall s\geq1,\,t\in\R/a\Z.
\end{equation}
\end{thm}
A consequence of (\ref{eq:d A u lam}) is that $|d_Au|_0$ decays as $e^{\left(-\frac{2\pi}a+\eps\right)s}$, for every $\eps>0$. Here the point-wise norm $|\cdot|_0$ is taken \wrt the standard metric $ds^2+dt^2$ on $\Si$ and the metric $g_{\om,J}$ on $M$. This generalizes a known decay result for pseudo-holomorphic maps, which corresponds to the case $\G=\{\one\}$ (see for example chapter 4 in \cite{MS}). In this case the result is optimal in the sense that $|du|_0$ does in general not decay faster than $e^{-\frac{2\pi}as}$. To see this consider $M:=\CP^1$ with the Fubini-Studi form, and let $u:\Si\to\CP^1\iso\C\cup\{\infty\}$ be defined by $u(z):=e^{\frac{2\pi}a z}$. On the other hand in some special examples the energy density $e_w$ decays faster than stated in Theorem \ref{thm:decay}, see for example the book \cite{JT} by Jaffe and Taubes. The overall strategy for the proof of Theorem \ref{thm:decay} is taken from the proof of Proposition 11.1 in \cite{GS}. The proof relies on an identity relating the energy of a vortex on a compact cylinder with the actions of its end-loops, and on the isoperimetric inequality of Theorem \ref{thm:isoperi}. The next result is an immediate consequence of Theorem \ref{thm:decay}, setting $a:=2\pi$ and applying the change of coordinates $\Si\ni z\mapsto e^{\frac{2\pi}az}\in \C.$ 
\begin{cor}[Optimal decay for vortices on the complex plane]\label{cor:decay C} Assume that hypothesis {\rm (H)} is satisfied, that $\Si=\C$, and that $\om_\Si$ is the standard area form on $\C$. Let $p>2$ and $w:=(u,A)$ be a locally $W^{1,p}$-solution of the equations {\rm (\ref{eq:vort P})} such that $E(w)<\infty$ and $u(P)\sub M$ has compact closure. Then for every $\eps>0$ there exists a constant $C$ such that 
   \begin{equation}\nn e_w(z)\leq C|z|^{-4+\eps},\quad\forall z\in\C\wo B_1.
  \end{equation}
 \end{cor}
I would like to thank the referees and Dietmar Salamon for some useful comments, Urs Frauenfelder and Kai Cieliebak for enlightening discussions, and Andreas Ott and Jan Swoboda for meticulous reading. This article was written while I was at the Ludwig-Maximilians University in Munich, Germany. During this time, I was supported by a grant from the Swiss National Science Foundation (SNF-Beitrags-Nr. PBEZ2-111606), which I gratefully acknowledge.
\section{Invariant symplectic action}\label{sec:action}
\subsection{Background}\label{subsec:back}
Let $(M,\om)$ be a symplectic manifold without boundary. We fix a Riemannian metric $\lan\cdot,\cdot\ran_M$ on $M$, and denote by $d,\exp,|v|,\iota_x>0$ and $\iota_X:=\inf_{x\in X}\iota_x\geq0$ the distance function, the exponential map, the norm of a vector $v\in TM$, and the injectivity radii of a point $x\in M$ and a subset $X\sub M$ respectively. We define the symplectic action of a loop $x:S^1\to M$ of length $\ell(x)<2\iota_{x(S^1)}$ to be
\begin{equation}
  \label{eq:A x int}\A(x):=-\int_\D u^*\om.
\end{equation}
Here $\D\sub\R^2$ denotes the (closed) unit disk, and $u:\D\to M$ is any smooth map such that 
\begin{equation}\label{eq:u e 2 pi}u(e^{2\pi it})=x(t),\,\forall t\in \R/\Z\iso S^1,\quad d\big(u(z),u(z')\big)<\iota_{x(S^1)},\,\forall z,z'\in\D.
\end{equation}
\begin{lemma}\label{le:A} The action $\A(x)$ is well-defined, \ie a map $u$ as above exists, and $\A(x)$ does not depend on the choice of $u$. 
\end{lemma}
\begin{proof} The lemma follows from an elementary argument, using the exponential map $\exp_{x(0+\Z)}:T_{x(0+\Z)}M\to M$. \end{proof}
Let now $M,\om,\G,\g,\mu$ and $\lan\cdot,\cdot\ran_M$ be as in Section \ref{sec:main}, and assume that $\lan\cdot,\cdot\ran_M$ is $\G$-invariant. The invariant action functional (also denoted by $\A$)
\[\A:\big\{(x,\xi)\in C^\infty(S^1,M\x\g)\,\big|\,x\textrm{ admissible}\big\}\to\R\]
is now defined as in Definition \ref{defi:action}. Note that if $x$ is admissible then $\A(x,\xi)$ is well-defined, \ie the required gauge transformation $g$ exists and the right hand side of (\ref{eq:A x xi}) does not depend on the choice of $g$. The functional $\A$ is invariant under the action of the gauge group $C^\infty(S^1,\G)$ on $C^\infty(S^1,M\x\g)$ given by $g_*(x,\xi):=\big(gx,(g\xi-\dot g)g^{-1}\big)$. To understand the definition of $\A(x,\xi)$ better, note that if both the lenghts $\ell(x)$ and $\ell(gx)$ are small then the term $-\int_0^1\big\lan\mu\circ x,g^{-1}\dot g\big\ran\,dt$ in (\ref{eq:A x xi}) compensates the effect of taking $\A(gx)$ rather than $\A(x)$. This is made precise in Lemma \ref{le:key} (\ref{le:key:A}) below. 

{\bf Example.} Consider $M:=\R^2$ endowed with the standard symplectic and Riemannian structures, and the action of the circle $\G:=S^1\iso\R/\Z$ by rotation, with moment map $\mu$ as in (\ref{eq:mu R 2}). Then every $x\in C^\infty\big(S^1,\R^2\wo\{0\}\big)$ is admissible, and the constant loop $x:\const0$ is inadmissible. If $|x|\const r>0$ and $\xi\const0$ then $A(x,\xi)=\pi(1-r^2)\deg\left(x/|x|:S^1\to S^1\right)$. 

The functional $\A$ is a modified version of a functional introduced by Gaio and Salamon in \cite{GS} (page 74). Their definition is based on the following lemma. Here for $\phi\in\g^*$ we denote by $|\phi|$ the norm of $\phi$ \wrt the inner product on $\g^*$ induced by $\lan\cdot,\cdot\ran_\g$. Recall that $L_x:\g\to T_xM$ denotes the infinitesimal action at $x\in M$. 
\begin{lemma}[Lemma 11.2 in \cite{GS}]\label{le:de c} Assume that $\mu$ is proper and that $\G$ acts freely on $\mu^{-1}(0)$. Then there are positive constants $\de$ and $c$ such that the following holds. If $x:S^1\to M$ and $\xi:S^1\to \g$ are smooth loops such that $\max_{S^1}|\mu\circ x|<\de$, then there is a point $x_0\in\mu^{-1}(0)$ and a smooth loop $g:S^1\to\G$ such that
\begin{equation}\label{eq:max d}\max_{S^1}\big|\xi+\dot gg^{-1}\big|\leq c\ell(x,\xi),\quad d(x(t),g(t)x_0)\leq c\big(|\mu\circ x(t)|+\ell(x,\xi)\big),
\end{equation}
where $\ell(x,\xi):=\int_0^1|\dot x+L_x\xi|\,dt$ denotes the twisted length of $(x,\xi)$.
\end{lemma}
Let $\de$ and $c$ be as in Lemma \ref{le:de c}, and suppose that $2c\de<\iota_M$. Assume that $(x,\xi)\in C^\infty(S^1,M)$ is a loop such that $\max_{S^1}|\mu\circ x|<\de$ and $\ell(x,\xi)<\de$. Then Gaio and Salamon define
\begin{equation}\label{eq:A 0 1}\A(x,\xi) := -\int_{[0,1]\x S^1}u^*\om+\int_0^1\big\lan\mu(x(t)),\xi(t)\big\ran\,dt.\end{equation}
Here $u:[0,1]\x S^1\to M$ is defined by $u(s,t):=\exp_{g(t)x_0}sv(t),$ where $x_0$ and $g$ are as in Lemma \ref{le:de c}, and $v(t)\in T_{g(t)x_0}M$ is the unique small tangent vector such that $x(t)=\exp_{g(t)x_0}v(t)$. Under the choices of $\de$ and $c$ above this expression may actually depend on the choice of $x_0$ and $g$, as an explicit example based on the $S^1$-action on $\R^2$ by rotation shows. One can overcome this difficulty by proving that $\de$ can be shrunk such that the following is satisfied. Namely, let $x$ and $\xi$ be as in the hypothesis of Lemma \ref{le:de c}, $(x_0,g)$ be as in the statement, $(x_0',g')$ be another such pair, and let $u$ and $u'$ be the corresponding maps defined as above. Then the triples $(x_0,g,u)$ and $(x_0',g',u')$ are smoothly homotopic. One should also replace the injectivity radius $\iota_M$ by $\iota_{\mu^{-1}(0)}$, since the former may be 0. Note that after these modifications the expression $\A(x,\xi)$ may depend on the pair $(\de,c)$ (assuming that it is well-defined).
\begin{rmk}\rm  Heuristically, the two expression (\ref{eq:A x xi}) and (\ref{eq:A 0 1}) for $\A(x,\xi)$ should be the same, provided that $\ell(x,\xi)+\max_{S^1}|\mu\circ x|$ is small enough. The idea how to see this is to choose $x_0$ and $g$ as in Lemma \ref{le:de c}. Then by the second inequality in (\ref{eq:max d}) the gauge transformed loop $x':=g^{-1}x$ is close to the point $x_0$. We define $\wt g:[0,1]\x S^1\to\G$ by $\wt g(s,t):=g(t)^{-1}$. Let $u$ be as in the definition (\ref{eq:A 0 1}) of $\A(x,\xi)$. Then the image of the map $u':=\wt gu:[0,1]\x S^1\to M$ lies inside a small ball around $x_0$. Since $u'$ maps the left part of the boundary of $[0,1]\x S^1$ to the point $x_0$, it induces a map from the disk $\D$ to $M$, whose restriction to the boundary circle agrees with $x'$. It follows that 
\begin{equation}\label{eq:A x'}\A(x')=-\int_{[0,1]\x S^1}{u'}^*\om=-\int_{[0,1]\x S^1}u^*\om+\int_{[0,1]\x S^1}d\big\lan\mu\circ u,\wt g^{-1}d\wt g\big\ran,\end{equation}
where in the second equality we used Lemma \ref{le:u' u} below. By Stokes' Theorem, the second term on the right hand side equals
\[\int_{S^1}\big\lan\mu\circ u,\wt g^{-1}d\wt g\big\ran\big|_{s=1}-\int_{S^1}\big\lan\mu\circ u,\wt g^{-1}d\wt g\big\ran\big|_{s=0}=\int_0^1\left\lan\mu\circ x,g\frac d{dt}(g^{-1})\right\ran\,dt - 0.\]
Here we used the fact that $u(0,t)=g(t)x_0\in\mu^{-1}(0)$. Combining this equality with (\ref{eq:A x'}), it follows that the right hand side of (\ref{eq:A x xi}) with $g$ replaced by $g^{-1}$ equals the right hand side of (\ref{eq:A 0 1}). The problem with this ``proof'' is that it is not clear how small $\ell(x,\xi)+\max_{S^1}|\mu\circ x|$ has to be in order for the arguments to work. 
\end{rmk}
\subsection{Key Lemma}\label{subsec:key}
Let $M,\om,\G,\mu,\lan\cdot,\cdot\ran_M,M^*,\ell$ and $\bar\ell$ be as in Section \ref{sec:main}. We denote by 
\begin{equation}\label{eq:Pr}\Pr:TM\to\big\{(x,L_x\xi)\,\big|\,x\in M,\,\xi\in\g\big\}\sub TM\end{equation}
the fiber-wise orthogonal projection to the image of the infinitesimal action.
\begin{lemma}[Key Lemma]\label{le:key} For every compact subset $K\sub M^*$ the following statements hold.
  \begin{enui}\item \label{le:key:C}  There exist constants $\de>0$ and $C$ with the following property. If $s_-\leq s_+$ are real numbers and $u\in C^\infty\big([s_-,s_+]\x S^1,K\big)$ is a map satisfying $\bar\ell(\G u(s,\cdot))<\de$ for every $s\in[s_-,s_+]$ then there exists a map $g\in C^\infty\big([s_-,s_+]\x S^1,\G\big)$ such that 
\begin{equation}
\label{eq:Pr d t}\big|\Pr\d_t(gu)(s,t)\big|\leq C\bar\ell(\G u(s,\cdot))^2,\quad \forall (s,t)\in[s_-,s_+]\x S^1.
\end{equation}
\item \label{le:key:A} There exists a number $0<\eps\leq 2\iota_K$ such that 
for every pair $(x,g)\in C^\infty(S^1,K\x\G)$ we have
\begin{equation}\label{eq:ell x eps}\ell(x)<\eps,\quad \ell(gx)<\eps \Then \A(gx)-\A(x)=\int_0^1\lan\mu\circ x,g^{-1}\dot g\ran\,dt.
\end{equation}
  \end{enui}
\end{lemma}
The following lemma is used in the proof of part (\ref{le:key:A}) of Lemma \ref{le:key}.
\begin{lemma}\label{le:u' u} Let $(M,\om)$ be a symplectic manifold, $\G$ a connected Lie group acting on $M$ in a Hamiltonian way, with moment map $\mu:M\to \g^*$, $X$ a manifold, and $u\in C^\infty(X,M)$ and $g\in C^\infty(X,\G)$ be maps. Then 
\[(gu)^*\om=u^*\om-d\big\lan\mu\circ u,g^{-1}dg\big\ran.\]
\end{lemma}
\begin{proof} Let $M,\om,G,\mu,X,u$ and $g$ be as in the hypothesis. Then 
\[\begin{array}{l}(gu)^*\om\\ 
=\om\big(gdu\cdot,gdu\cdot\big)+\om\big((dg\cdot)u,(dg\cdot)u\big)+\om\big((dg\cdot)u,gdu\cdot\big)+\om\big(gdu\cdot,(dg\cdot)u\big)\\
=u^*\om+\om\big(L_ug^{-1}dg\cdot,L_ug^{-1}dg\cdot\big)+\om\big(L_ug^{-1}dg\cdot,du\cdot\big)+\om\big(du\cdot,L_ug^{-1}dg\cdot\big)\\
=u^*\om+\frac12\big\lan\mu\circ u,\big[g^{-1}dg\wedge g^{-1}dg\big]\big\ran-\big\lan d\mu(u)du\wedge g^{-1}dg\big\ran\\
=u^*\om-d\big\lan\mu\circ u,g^{-1}dg\big\ran.
\end{array}\]
Here we used the notation $\eta x :=L_{gx}(\eta g^{-1})$, $gv:=\left.\frac d{dt}\right|_{t=0}g\ga(t)\in T_{gx}M$, for $x\in M$, $g\in\G$, $\eta\in T_g\G$ and $v\in T_xM$, where $\ga\in C^\infty(\R,M)$ is a curve satisfying $\ga(0)=x$, $\dot\ga(0)=v$. Furthermore, $\om\big((dg\cdot)u,(dg\cdot)u\big)$ denotes the two-form $TX\x TX\ni(\ze,\ze')\mapsto\om\big((dg\,\ze)u,(dg\,\ze')u\big)$, and similarly for the other expressions. This proves Lemma \ref{le:u' u}. 
\end{proof}
\begin{proof}[Proof of Lemma \ref{le:key}]\setcounter{claim}{0} Let $K\sub M^*$ be a compact subset. We may assume \Wlog that $M=M^*$ and $K$ is $\G$-invariant. We fix an ad-invariant inner product $\lan\cdot,\cdot\ran_\g$ on $\g$, and denote by $d^\G$ and $\iota^\G$ the corresponding distance function on $\G$ and injectivity radius of $\G$.

{\bf We prove (\ref{le:key:C}).} Consider the compact subset in the quotient $\bar K:=K/\G\sub M/\G.$  We choose an open neighborhood $X\sub M/\G$ of $\bar K$ with compact closure, and we equip it with the Riemannian metric induced by $\lan\cdot,\cdot\ran_M$. We denote by $P\sub M$ the pre-image of $X$ under the canonical projection $M\to M/\G$. This is a principal $\G$-bundle with right-action given by $P\x\G\to P,$ $(x,g)\mapsto g^{-1}x$. Applying Proposition \ref{prop:holon} of appendix \ref{sec:holon} with $K$ replaced by $\bar K\sub X$ there exists a constant $C$ such that the conclusion of this proposition holds. We define the connection one-form $A$ on $P$ by
\begin{equation}
  \label{eq:A}A_xv:=\left(L_x^*L_x\right)^{-1}L_x^*v\in \g,
\end{equation}
for $x\in P$ and $v\in T_xP$. Here $L_x^*:T_xP\to \g$ denotes the adjoint map to $L_x$ \wrt the metric on $P$ and the inner product $\lan\cdot,\cdot\ran$ on $\g$. It follows from the assumption $M=M^*$ that $L_x$ is injective, hence $A$ is well-defined. We choose $\de>0$ less than the injectivity radius of the subset $\bar K$ in $X$ and such that $C||F_A||_{L^\infty(X)}\de^2<\iota^\G$. Let $s_-,s_+$ and $u$ be as in the hypothesis of part (\ref{le:key:C}). We choose a smooth $t$-horizontal lift $\wt u:[s_-,s_+]\x\R\to K$ of the map $\G u:[s_-,s_+]\x S^1\to \bar K$. This means that $\G\wt u(s,t)=\G u(s,t+\Z)$ and $A\d_t\wt u\const0$. We fix $s\in [s_-,s_+]$, and define $\wt x:=\wt x_s:\R\to K$ by $\wt x(t):=\wt u(s,t)$, and $h:=h_s\in\G$ by the equation $\wt x(1)=:h\wt x(0)$. It follows that $h$ is the holonomy of $A$ around the loop $\G u(s,\cdot)$ with base point $\wt x(0)$. Thus by the assertion of Proposition \ref{prop:holon} and the inequality $C||F_A||_{L^\infty(X)}\de^2<\iota^\G$ we have
\begin{equation}\label{eq:d G C}d^\G(\one,h)\leq C||F_A||_{L^\infty(X)}\bar\ell(\G u(s,\cdot))^2<\iota^\G.
\end{equation}
Hence there exists a unique element $\xi:=\xi_s\in\g$ satisfying
\begin{equation}\label{eq:exp d}\exp\xi=h,\qquad |\xi|=d^\G(\one,h).
\end{equation}
We define the map $\wt h:=\wt h_s:\R\to\G$ by $\wt h(t):=\exp(-t\xi)$, and the map $\wt g:=\wt g_s:\R\to\G$ by the equation 
\begin{equation}\label{eq:u s t Z}(\wt g\wt h\wt x)(t)=u(s,t+\Z).
\end{equation}
\begin{claim}\label{claim:wt g} We have $\wt g(t+1)=\wt g(t)$. 
\end{claim}
\noindent{\bf Proof of Claim \ref{claim:wt g}.} By (\ref{eq:u s t Z}) and the equation $\wt x(1)=h\wt x(0)$ we have $\wt g(0)=\wt g(1)$. We show that the maps $\wt g$ and $t\mapsto \wt g(t+1)$ satisfy the same ordinary differential equation: Using (\ref{eq:u s t Z}) again, we get 
\begin{equation}\label{eq:A d}A\d_tu(s,\cdot+\Z)=A\Big(\dot{\wt g}\wt h\wt x-\wt g\xi\wt h\wt x+\wt g\wt h\dot{\wt x}\Big)=\dot{\wt g}\wt g^{-1}-\wt g\xi \wt g^{-1}+0.
\end{equation}
Here we used the notation of the proof of Lemma \ref{le:u' u}. Furthermore, in the second step we used the fact $A\d_t\wt u=0$. A similar calculation shows that $A\d_tu(s,\cdot+\Z)=\big(\dot{\wt g}\wt g^{-1}-\wt g\xi \wt g^{-1}\big)(\cdot+1)$. Combining this with (\ref{eq:A d}) we obtain Claim \ref{claim:wt g}.

Using Claim \ref{claim:wt g} we may define the map $g:[s_-,s_+]\x S^1\to\G$ by $g(s,t+\Z):=\wt g_s(t)^{-1}$. Fixing $s\in[s_-,s_+]$, equation (\ref{eq:u s t Z}) and the facts $\Pr_x=L_xA_x$ (for $x\in P$) and $A\dot{\wt x}=0$ imply that
\[\Pr\d_t(gu)(s,\cdot+\Z)=\Pr\big(-\xi\wt h\wt x+\wt h\dot{\wt x}\big)=-L_{(gu)(s,\cdot+\Z)}\xi.\]
Defining $C':=\max\left\{|L_x\eta|\,|\,x\in \BAR P,\,\eta\in\g:\,|\eta|\leq1 \right\}$, we get 
\[\big|\Pr\d_t(gu)(s,t+\Z)\big|\leq C'|\xi_s|= C'd^\G(\one,h_s)\leq C'C||F_A||_{L^\infty(X)}\bar\ell(\G u(s,\cdot))^2.\]
Here in the second step we used the second identity in (\ref{eq:exp d}), and in the third step we used (\ref{eq:d G C}). This completes the proof of statement (\ref{le:key:C}). 

{\bf We prove (\ref{le:key:A}).} We define 
\begin{eqnarray}\nn C_1&:=&\max\left\{|L_x\xi|\,\big|\,x\in M:\,d(x,K)\leq\iota_K/4,\,\xi\in\g:\,|\xi|\leq1\right\},\\
C_2&:=&\sup \left\{\frac{d^\G(\one,g)}{d(x,gx)}\,\bigg|\,\one\neq g\in\G,\,x\in K\right\},
\end{eqnarray}
and we choose a positive number $\eps$ satisfying
\begin{equation}\label{eq:eps min}\eps<\min\left\{\frac{\iota_K}2,\frac{\iota_K}{4C_1C_2},\frac{\iota^\G}{C_2}\right\}.
\end{equation}
Since $\G$ and $K$ are compact and by assumption $\G$ acts freely on $M$, it follows that $C_2<\infty$. Let $(x,g)\in C^\infty(S^1,K\x\G)$ be a pair of loops satisfying $\ell(x)<\eps$ and $\ell(gx)<\eps$. By replacing $x$ and $g$ by $g(0+\Z)x$ and $g\cdot g(0+\Z)^{-1}$, we may assume \Wlog that $g(0+\Z)=\one$. We fix a point $t\in\R/\Z$. Then
\begin{eqnarray}\nn d^\G(\one,g(t))&\leq&C_2d\big(x(t),(gx)(t)\big)\\
\nn&\leq&C_2\Big(d\big(x(t),x(0+\Z)\big)+d\big((gx)(0+\Z),(gx)(t)\big)\Big)\\
\nn&\leq&\frac{C_2}2(\ell(x)+\ell(gx))\\
\label{eq:d G}&<&C_2\eps<\iota^\G.
\end{eqnarray}
Hence there exists a unique element $\xi_0(t)\in B_{\iota^\G}\sub \g$ such that 
\begin{equation}
  \label{eq:xi 0}\exp\xi_0(t)=g(t),\qquad |\xi_0(t)|=d^\G(\one,g(t)).
\end{equation}
Since $\ell(x)<\eps<\frac{\iota_K}2$, there exists a smooth map $u:\D\to B_{\iota_K/4}(x(0+\Z))\sub M$ such that $u(e^{2\pi it})=x(t)$ for every $t\in\R/\Z$. We choose a smooth function $\rho:[0,1]\to[0,1]$ that vanishes in a neighborhood of 0, and equals 1 in a neighborhood of 1, and we define $h:\D\to \G$ by $h(re^{2\pi it}):=\exp\big(\rho(r)\xi_0(t)\big)$. Furthermore, we define $u':\D\to M$ by $u'(z):=h(z)u(z)$.
\begin{claim}\label{claim:u'} We have $\A(gx)=-\int_\D {u'}^*\om$.
\end{claim}
\noindent{\bf Proof of Claim \ref{claim:u'}.} Using (\ref{eq:xi 0}) and $u(e^{2\pi it})=x(t)$, we have $u'(e^{2\pi it})=g(t)x(t)$, for $t\in\R/\Z$. Furthermore, we fix $r\in[0,1]$ and $t\in\R$, and define the path $\ga:[0,1]\to M$ by $\ga(\lam):=\exp\big(\lam \rho(r)\xi_0(t)\big)u(re^{2\pi it})$. Then
\[d\big(u(re^{2\pi it}),u'(re^{2\pi it})\big)\leq \ell(\ga)\leq C_1|\xi_0(t)|= C_1d^\G(\one,g(t))<C_1C_2\eps\leq\iota_K/4,\]
where in the third step we used (\ref{eq:xi 0}), in the fourth step (\ref{eq:d G}), and in the last step (\ref{eq:eps min}). Since $u'(1)=u(1)$, it follows that 
\[d(u'(1),u'(z))\leq d(u(1),u(z))+d\big(u(z),u'(z)\big)<\iota_K/2,\]
for $z\in\D$. Hence $d(u'(z),u'(z'))<\iota_K$, for $z,z'\in\D$. Claim \ref{claim:u'} follows now from the definition of $\A(gx)$.

By Claim \ref{claim:u'} and Lemma \ref{le:u' u}, we get
\[\A(gx)=\int_\D\Big(-u^*\om+d\big\lan\mu\circ u,h^{-1}dh\big\ran\Big)=\A(x)+\int_{\d\D}\big\lan\mu\circ u,h^{-1}dh\big\ran.\]
Here the second step follows from (\ref{eq:A x int}), which holds since $u$ satisfies (\ref{eq:u e 2 pi}). The equality stated in (\ref{eq:ell x eps}) follows now from the fact $h(e^{2\pi it})=g(t)$. This proves part (\ref{le:key:A}) and completes the proof of Lemma \ref{le:key}.
\end{proof}
\subsection{Proof of the isoperimetric inequality}\label{sec:proof isoperi}
Given a pair of loops $(x,\xi)\in C^\infty(S^1,K\x\g)$, the idea of the proof of Theorem \ref{thm:isoperi} is to gauge transform $x$ to a short loop $x'$. This is possible by the Key Lemma, under the assumption that the loop $\G x:S^1\to M^*/\G$ is short. The equivariant isoperimetric inequality then follows from the isoperimetric inequality for the case $\G=\{\one\}$.

\begin{proof}[Proof of Theorem \ref{thm:isoperi}]\setcounter{claim}{0}\label{isoperi:proof} Let $J$ be a $\G$-invariant $\om$-compatible almost complex structure on $M$ such that $g_{\om,J}=\lan\cdot,\cdot\ran_M$. The subset $M^*\sub M$ of all points on which $\G$ acts freely is open. Hence we may assume \Wlog that the action of $\G$ on $M$ is free. Let $K\sub M$ be a compact subset and $c>\frac1{4\pi}$ be a constant. By replacing $K$ by the compact set $\G K$ we may assume that $K$ is $\G$-invariant. We choose a constant $c_0\in\big(\frac1{4\pi},c\big).$ Applying Theorem 4.4.1 of \cite{MS} (isoperimetric inequality for the usual action) with $c$ replaced by $c_0$, there exists a constant $0<\de_0<\iota_K$ such that 
  \begin{equation}
    \label{eq:A x c 0}
\A(x)\leq c_0\ell(x)^2,
  \end{equation}
for every loop $x\in C^\infty(S^1,K)$ of length $\ell(x)<\de_0$. (Strictly speaking, in Theorem 4.4.1 in \cite{MS} it is assumed that $M$ is compact. However, the proof of this theorem carries over to the present situation.) Moreover, let $\de_1$ and $C_1$ be constants as in Lemma \ref{le:key}(\ref{le:key:C}) corresponding to $\de$ and $C$, and let $\eps$ be as in part (\ref{le:key:A}) of that lemma. We choose $\de>0$ such that
\begin{equation}\label{eq:1 de}\displaystyle\de<\min\left\{\frac{\de_0}2,\de_1,\frac\eps2,\frac1{C_1}\right\},\quad\displaystyle\Big(\sqrt{1+C_1\de}+\sqrt{(1+2C_1\de)C_1\de}\Big)^2 <\frac c{c_0}.\end{equation}
Let $x\in C^\infty(S^1,K)$ be a loop such that $\bar\ell(\G x)<\de$. By the assertion of Lemma \ref{le:key}(\ref{le:key:C}) there is a loop $g\in C^\infty(S^1,\G)$ such that 
\begin{equation}\label{eq:P d dt}\left|\left|\Pr\frac d{dt}(gx)\right|\right|_\infty\leq C_1\bar\ell(\G x)^2.
\end{equation}
We define $x':=gx$. It follows that 
\begin{equation}\label{eq:ell x'}\ell(x')\leq\big|\big|(\id-\Pr)\dot x'\big|\big|_1+ ||\Pr \dot x'||_1\leq\bar\ell(\G x)+C_1\bar\ell(\G x)^2<\de+\de,\end{equation}
where in the third step we used the fact that $\de\leq C_1^{-1}$.
\begin{claim}\label{claim:adm} The loop $x$ is admissible.
\end{claim}
\noindent{\bf Proof of Claim \ref{claim:adm}.} By (\ref{eq:ell x'}) and the inequalities $2\de<\eps\leq 2\iota_K\leq2\iota_{x(S^1)}$ the condition $\ell(x')<2\iota_{x(S^1)}$ is satisfied. Let now $\wt g\in C^\infty(S^1,\G)$ be a loop such that, setting $\wt x:=\wt gx$, we have $\ell(\wt x)\leq\ell(x').$ Applying Lemma \ref{le:key}(\ref{le:key:A}) with $x,g$ replaced by $x',\wt gg^{-1}$, we get
\begin{eqnarray}\nn \A(\wt x)-\A(x')=\int_0^1\left\lan\mu\circ x',g\wt g^{-1}\frac d{dt}(\wt g g^{-1})\right\ran\,dt=\int_0^1\big\lan\mu\circ x,\wt g^{-1}\dot{\wt g}-g^{-1}\dot g\big\ran\,dt.
\end{eqnarray}
Hence condition (\ref{eq:A wt g x int}) is satisfied. This proves Claim \ref{claim:adm}.

Let $\xi\in C^\infty(S^1,\g)$ be a loop, and $p\in[1,2]$.  
\begin{claim}\label{claim:isoperi} The isoperimetric inequality (\ref{eq:isoperi}) holds. 
\end{claim}
\noindent{\bf Proof of Claim \ref{claim:isoperi}.} We define $\xi':=(g\xi-\dot g)g^{-1}$ and $p':=\frac p{p-1}\in[2,\infty]$. By (\ref{eq:ell x'}) and (\ref{eq:1 de}) we have $\ell(x')<\de_0$, hence by (\ref{eq:A x c 0}) with $x$ replaced by $x'$,
\begin{eqnarray}\nn |\A(x',\xi')|&\leq&|\A(x')|+\left|\int_0^1\lan\mu\circ x',\xi'\ran\,dt\right|\\
\nn&\leq&c_0\ell(x')^2+||\xi'||_p\,||\mu\circ x'||_{p'}\\
\nn&\leq&c_0||\dot x'||_p^2+c_0m_K^2||\xi'||_p^2+\frac{1}{4c_0m_K^2}||\mu\circ x'||_{p'}^2\\
\label{eq:A c 0}&\leq&c_0\big(||\dot x'||_p^2+||L_{x'}\xi'||_p^2\big)+\frac\pi{m_K^2}||\mu\circ x'||_{p'}^2\,.
\end{eqnarray}
Here  the constant $m_K$ appearing in the third step is defined as in (\ref{eq:m X}) with $X:=K$, and in the fourth step we used the definition of $m_K$ and the fact that $c_0>\frac1{4\pi}$. 
\begin{claim}\label{claim:dot x '} The inequality $||\dot x'||_p^2+||L_{x'}\xi'||_p^2\leq\frac c{c_0}\big|\big|\dot x'+L_{x'}\xi'\big|\big|_p^2$ holds. 
\end{claim}
\noindent{\bf Proof of Claim \ref{claim:dot x '}.} Since $p\leq 2$, the map $||\cdot||:\R^2\to \R$, $||v||:=\left(v_1^{2/p}+v_2^{2/p}\right)^{p/2}$ is a norm. Hence, defining $f:=\big(|\dot x'|^p,|L_{x'}\xi'|^p\big):S^1\to \R^2$, we obtain 
\begin{eqnarray}\nn||\dot x'||_p^2+||L_{x'}\xi'||_p^2&=&\left|\left|\int_{S^1}f\,dt\right|\right|^{\frac2p}\\
\nn&\leq&\left(\int_{S^1}||f||\,dt\right)^{\frac2p}\\
\nn&\leq& \left|\left|\sqrt{|\dot x'|^2+|L_{x'}\xi'|^2}\right|\right|_p^2\\ 
\label{eq:dot x '}&=&\left|\left|\sqrt{\big|\dot x'+L_{x'}\xi'\big|^2-2g_{\om,J}\big(\dot x',L_{x'}\xi'\big)}\right|\right|_p^2.
\end{eqnarray} 
Furthermore, we have
\begin{eqnarray}\nn \big|g_{\om,J}\big(\dot x',L_{x'}\xi'\big)\big|&=&\big|g_{\om,J}\big( \Pr\dot x',-\Pr\dot x'+\dot x'+L_{x'}\xi'\big)\big|\\
\nn &\leq&|\Pr\dot x'|^2+\big|g_{\om,J}\big(\Pr\dot x',\dot x'+L_{x'}\xi'\big)\big|\\
\label{eq:x' L} &\leq&\left(1+\frac1{2C_1\de}\right)|\Pr\dot x'|^2+\frac{C_1\de}2|\dot x'+L_{x'}\xi'|^2,
\end{eqnarray}
where in the last step we used Young's inequality. Moreover, inequality (\ref{eq:P d dt}) and the fact $\bar\ell(\G x)<\de$ imply that
\begin{equation}\label{eq:Pr dot x '}|\Pr\dot x'|\leq C_1\bar\ell(\G x)^2\leq C_1\de\big|\big|\dot x'+L_{x'}\xi'\big|\big|_1\leq C_1\de\big|\big|\dot x'+L_{x'}\xi'\big|\big|_p\,.
\end{equation}
We define $c_1:=\sqrt{1+C_1\de}$ and $c_2:=\sqrt{(1+2C_1\de)C_1\de}$. Combining (\ref{eq:x' L}) and (\ref{eq:Pr dot x '}), we obtain
\[\big|\dot x'+L_{x'}\xi'\big|^2-2g_{\om,J}\big(\dot x',L_{x'}\xi'\big)\leq c_1^2 \big|\dot x'+L_{x'}\xi'\big|^2+c_2^2\big|\big|\dot x'+L_{x'}\xi'\big|\big|_p^2.\]
Combining this with (\ref{eq:dot x '}), we get 
\begin{eqnarray}\nn ||\dot x'||_p^2+||L_{x'}\xi'||_p^2&\leq&\left|\left|c_1 \big|\dot x'+L_{x'}\xi'\big|+c_2\big|\big|\dot x'+L_{x'}\xi'\big|\big|_p\right|\right|_p^2\\
\nn&\leq&\big((c_1+c_2)\big|\big|\dot x'+L_{x'}\xi'\big|\big|_p\big)^2\\
\nn&\leq&\frac c{c_0}\big|\big|\dot x'+L_{x'}\xi'\big|\big|_p^2.
\end{eqnarray}
Here in the second step we used Minkowski's inequality and the fact that the Haar measure of $S^1$ is 1. Furthermore, in the last step we used (\ref{eq:1 de}). This proves Claim \ref{claim:dot x '}. 

Claim \ref{claim:isoperi} follows from (\ref{eq:A c 0}), Claim \ref{claim:dot x '} and the equalities 
\[\A(x,\xi)=\A(x',\xi'),\quad \big|\dot x'+L_{x'}\xi'\big|=\big|\dot x+L_x\xi\big|,\quad |\mu\circ x'|=|\mu\circ x|.\] 
This completes the proof of Theorem \ref{thm:isoperi}.\end{proof}
\begin{rmk}\rm \label{rmk:alt} There is an alternative approach to the isoperimetric inequality for $p=2$. Namely, as pointed out to me by Urs Frauenfelder, we may interpret the invariant action as a Morse-Bott function $f$, defined on the infinite dimensional space $X$ of gauge equivalence classes of loops $(x,\xi)\in C^\infty(S^1,M^*\x\g)$, for which $x$ is admissible. Assuming that hypothesis (H) above is satisfied, the set $\Crit f$ of critical points of $f$ can be identified with the symplectic quotient via the map
\[\mu^{-1}(0)/\G\ni \G x\mapsto [x,0]\in\Crit f,\]
where $[x,0]$ denotes the equivalence class of the constant map $(x,0)$. Since $\Crit f\sub f^{-1}(0)$, heuristically, for every constant 
\[ c>\big(2\min\big\{|\lam|\,\big|\,\lam\textrm{ eigenvalue of }H_f^\perp(p),\,p\in\Crit f\big\}\big)^{-1}\]
there exists a neighborhood $U$ of $\Crit f$ such that for every $p\in U$ we have
\begin{equation}\label{eq:f p}|f(p)|\leq c|\na f(p)|^2.\end{equation}
Here we fix a Riemannian metric $\lan\cdot,\cdot\ran_X$ on $X$, and denote by $|\na f(p)|$ the corresponding norm of the gradient of $f$ at a point $p=[x,\xi]\in X$, and by $H_f^\perp(p)$ the Hessian of $f$ at $p$. Choosing a suitable metric $\lan\cdot,\cdot\ran_X$, the isoperimetric inequality (\ref{eq:isoperi}) with $p:=2$ can be derived from (\ref{eq:f p}). 
\end{rmk}
\section{Symplectic vortices}\label{sec:vort}
\subsection{Energy action identity and bound on energy density}\label{subsec:en act bound}
Let $M,\om,$ $\G,\g,\lan\cdot,\cdot\ran_\g,\mu$ and $J$ be as in Section \ref{sec:main}, $a>0$ be a number, and $\Si:=\big\{s+it\in \C\,\big|\,s\geq0\big\}/ai\Z$ be the half-cylinder, equipped with an $i$-compatible area form $\om_\Si=\lam^2ds\wedge dt$. In this and the next subsection we identify $\g^*$ with $\g$ via $\lan\cdot,\cdot\ran_\g$. Furthermore, all norms of vectors in $TM$ \etc are \wrt the metric $g_{\om,J}$. Let $P\to\Si$ be a principal $\G$-bundle. Since by assumption $\G$ is connected the bundle $P$ is trivial. Hence for the proof of Theorem \ref{thm:decay} it suffices to consider the case $P=\Si\x\G$. In this case, an equivariant map from $P$ to $M$ corresponds in a bijective way to a map $u:\Si\to M$. Furthermore, a connection one-form $A$ on $P$ bijectively corresponds to a pair of maps $\Phi,\Psi:\Si\to\g,$ via the formula
\begin{equation}
  \label{eq:A z g}A_{(z,g)}(\ze_1,\ze_2,g\xi)=\xi+g^{-1}\big(\ze_1\Phi(z)+\ze_2\Psi(z)\big)g,
\end{equation}
for $(z,g)\in P$ and $(\ze_1+i\ze_2,g\xi)\in T_{(z,g)}P=\C\x g\cdot\g$. The equations (\ref{eq:vort P}) are equivalent to the $\lam$-vortex  equations
 \begin{equation}
  \label{eq:vort lam}\left\{
\begin{array}{ccc}\d_su+L_u\Phi+J\big(\d_tu+L_u\Psi\big)&=&0\\
\big(\d_s\Psi-\d_t\Phi+[\Phi,\Psi]\big)+\lam^2\mu\circ u&=&0.
  \end{array}\right.
\end{equation}
For simplicity, we restrict now to the case $a=1$, and we identify
\[S^1\iso\R/\Z,\qquad \Si=\big\{s+it\in \C\,\big|\,s\geq0\big\}/i\Z\iso\R\x S^1.\]Given an open subset $U\sub \Si$ and a number $p>2$ the energy density and the energy of a solution $w:=(u,\Phi,\Psi)\in W^{1,p}_\loc(U,M\x\g\x\g)$ of the equations (\ref{eq:vort lam}), \wrt the standard metric $ds^2+dt^2$ on $\Si$, are given by
\begin{equation}
  \label{eq:wt e w} \wt e_w=|\d_su+L_u\Phi|^2+\lam^2|\mu\circ u|^2,\qquad E(w,U)=\int_U\wt e_w\,ds\wedge dt.
\end{equation}
\begin{prop}[Energy action identity]\label{prop:en act} For every compact subset $K\sub M^*$ there exists a constant $\de>0$ with the following property. Let $s_-\leq s_+$ be numbers, $\Si:=[s_-,s_+]\x S^1$ be the compact cylinder, $\lam\in C^\infty\big(\Si,(0,\infty)\big)$ be a function, and $w:=(u,\Phi,\Psi)\in C^\infty\big(\Si,K\x\g\x\g\big)$ be a solution of the equations {\rm (\ref{eq:vort lam})} satisfying $\bar\ell(\G u(s,\cdot))<\de$ for every $s\in [s_-,s_+]$. Then the loops $u(s_-,\cdot)$ and $u(s_+,\cdot)$ are admissible, and
  \begin{equation}
\label{eq:E A A} E(w,\Si)=-\A\big((u,\Psi)(s_+,\cdot)\big)+\A\big((u,\Psi)(s_-,\cdot)\big).
  \end{equation}
\end{prop}
For the proof of Proposition \ref{prop:en act} we need the following lemma. Recall that by $\iota_x>0$ we denote the injectivity radius of a point $x\in M$.
\begin{lemma}\label{le:en act} Let $s_-\leq s_+$ be numbers, $\Si:=[s_-,s_+]\x S^1$ be the compact cylinder, $\lam\in C^\infty\big(\Si,(0,\infty)\big)$ be a function, and let $w:=(u,\Phi,\Psi)\in C^\infty\big(\Si,M\x\g\x\g\big)$ be a solution of the equations {\rm (\ref{eq:vort lam})}, satisfying 
  \begin{equation}
    \label{eq:ell u s}\ell(u(s,\cdot))<2\inf_{t\in S^1}\iota_{u(s,t)},\quad \forall s\in [s_-,s_+].
  \end{equation}
Then
\begin{eqnarray}\nn E\big(w,\Si\big)&=&-\A(u(s_+,\cdot))+\A((u(s_-,\cdot))\\ 
\label{eq:E A int u}&&+\int_0^1\Big(-\big\lan\mu\circ u,\Psi\big\ran\big|_{(s_+,t)}+\big\lan\mu\circ u,\Psi\big\ran\big|_{(s_-,t)}\Big)\,dt.
\end{eqnarray}
\end{lemma}
\begin{proof}[Proof of Lemma \ref{le:en act}]\setcounter{claim}{0} As in the proof of Proposition 3.1 in \cite{CGS} we have
\[\wt e_w=\om(\d_su,\d_tu)-\d_s\lan\mu\circ u,\Psi\ran+\d_t\lan\mu\circ u,\Phi\ran.\]
The lemma follows by integrating this over $\Si$ and using the energy action identity for the usual action. 
\end{proof}
\begin{proof}[Proof of Proposition \ref{prop:en act}]\setcounter{claim}{0} Let $K\sub M^*$ be a compact subset. We fix constants $\de,C,\eps$ as in Lemma \ref{le:key}, such that 
\begin{equation}\label{eq:de C}C\de^2+\de<\eps<2\iota_K=2\inf_{x\in K}\iota_x.
\end{equation}
Let $s_-\leq s_+,\lam$ and $w:=(u,\Phi,\Psi)$ be as in the hypothesis of Proposition \ref{prop:en act}. By the assertion of Lemma \ref{le:key}(\ref{le:key:C}) there exists a map $g\in C^\infty\big([s_-,s_+]\x S^1,\G\big)$ such that inequality (\ref{eq:Pr d t}) holds, where the projection $\Pr$ is defined as in (\ref{eq:Pr}). Hence fixing $s\in[s_-,s_+]$, we may estimate
\begin{eqnarray}\nn \ell\big((gu)(s,\cdot)\big)&\leq&\int_0^1\big|\Pr\d_t(gu)(s,t)\big|\,dt+\int_0^1\big|(\id-\Pr)\d_t(gu)(s,t)\big|\,dt\\
\nn&\leq&C\bar\ell(\G u(s,\cdot))^2+\bar\ell(\G u(s,\cdot))\\
\label{eq:ell g u s}&\leq&C\de^2+\de.
\end{eqnarray}
Combining this with (\ref{eq:de C}) we get $\ell\big((gu)(s,\cdot)\big)<2\iota_K\leq2\inf_{t\in S^1}\iota_{u(s,t)}$. Furthermore, let $\wt g\in C^\infty(S^1,\G)$ be a loop such that $\ell\big(\wt g u(s,\cdot)\big)\leq\ell\big((gu)(s,\cdot)\big)$. Since $\ell(u(s,\cdot))<\de<\eps$, Lemma \ref{le:key}(\ref{le:key:A}) implies that equality (\ref{eq:A wt g x int}) holds with $x$ and $g$ replaced by $u(s,\cdot)$ and $g(s,\cdot)$. It follows that $u(s,\cdot)$ is admissible. The equality (\ref{eq:E A A}) follows now from Lemma \ref{le:en act} with $w$ replaced by the gauge transformed map $g_*w$. This proves Proposition \ref{prop:en act}.
\end{proof}
\begin{lemma}[Point-wise bound on $\wt e_w$]\label{le:e w} Assume that hypothesis {\rm (H)} of Section \ref{sec:main} holds. Let $\Si=\big\{s+it\in \C\,\big|\,s\geq0\big\}/i\Z$ be the half-cylinder, $\om_\Si=\lam^2 ds\wedge dt$ be an area form on $\Si$ that satisfies the second inequality in (\ref{eq:lam sup}), and let $w:=(u,\Phi,\Psi)$ be a smooth solution of the equations {\rm (\ref{eq:vort lam})} on $\Si$ of finite energy $E(w,\Si)$, such that $\BAR{u(\Si)}\sub M$ is compact. Then there exists a number $s_0\geq\frac12$ such that for $z\in \big([s_0,\infty)+i\R\big)/i\Z$ we have
  \begin{equation}
\nn \wt e_w(z)\leq \frac{32}\pi E\big(w,B_{\frac12}(z)\big).
  \end{equation} 
\end{lemma}
The proof of Lemma \ref{le:e w} uses the following two lemmas. For $r>0$ we denote by $B_r\sub\R^2$ the open ball of radius $r$ around 0. 
\begin{lemma}[Mean value]\label{le:mean} Let $r>0$, $C\geq0$ and $f\in C^2(B_r,\R)$.
\[ \textrm{If}\quad f\geq0,\,\La f\geq-Cf^2,\, \int_{B_r}f<\frac\pi{8C}\Then f(0)\leq\frac8{\pi r^2}\int_{B_r}f.\]
\end{lemma}
\begin{proof}[Proof of Lemma \ref{le:mean}] This is Lemma 4.3.2. in \cite{MS}.
\end{proof}
Recall the definition (\ref{eq:m X}) of $m_X$. For $\de>0$ we denote by $\bar B_\de\sub\g$ the closed ball of radius $\de$ around the origin. 
\begin{lemma}\label{le:c c} Assume that hypothesis {\rm (H)} holds. Then for every $c<m_{\mu^{-1}(0)}$ there exists a number $\de>0$ such that $\G$ acts freely on $\mu^{-1}(\bar B_\de)$, and 
\[|L_x\xi|\geq c|\xi|,\quad \forall x\in\mu^{-1}(\bar B_\de),\, \xi\in\g.\] 
\end{lemma}
\begin{proof} This follows from an elementary argument involving the sequence $K_\nu:=\mu^{-1}\left(\bar B_{1/\nu}\right)\sub M$, for $\nu\in\N$.
\end{proof}
\begin{proof}[Proof of Lemma \ref{le:e w}] Let $\Si$, $\om_\Si=\lam^2ds\wedge dt$ and $w$ be as in the hypothesis.
  \begin{claim}\label{claim:La e w lam}There exists a constant $C$ such that $\La \wt e_w\geq-C\wt e_w^2.$
\end{claim}
\noindent{\bf Proof of Claim \ref{claim:La e w lam}.} We abbreviate $v_s:=\d_su+L_u\Phi$ and $v_t:=\d_tu+L_u\Psi$, and define
\begin{eqnarray}\nn I&:=&2\lam^{4}|L_u\mu\circ u|^2+4\lam^{2}|L_u^*v_s|^2+4\lam^{2}|L_u^*v_t|^2\\
\nn II&:=&6\d_t(\lam^{2})\lan\mu\circ u,L_u^*v_s\ran-6\d_s(\lam^{2})\lan\mu\circ u,L_u^*v_t\ran+\La(\lam^{2})|\mu\circ u|^2.
\end{eqnarray}
Let $a>0$. Then it follows from formula (90) on page 63 in the article \cite{GS} by R. Gaio and D. A. Salamon, with $\eps:=1$, that there exists a constant $b>0$ (depending on the compact set $\BAR{u(\Si)}\sub M$) such that 
\begin{equation}
  \label{eq:III}III:=\La \wt e_w-I-II\geq -a\lam^4\big|\mu\circ u\big|^2-b|v_s|^4.
\end{equation}
We define $C:=\sup_{\Si}\big(|d(\lam^{-1})|^2+\La(\lam^{-2})\big)$. Then by the second condition in (\ref{eq:lam sup}) we may choose a constant $c<m_{\mu^{-1}(0)}$ such that $C<2c^2$. We set $a:=2c^2-C$ and choose a constant $b$ as above. By Young's inequality with exponent 2, we have
\begin{eqnarray}\nn 6\d_t(\lam^2)\lan\mu\circ u,L_u^*v_s\ran&\geq&-\frac{9\big(\d_t(\lam^2)|\mu\circ u|\big)^2}{4\lam^2}-4\lam^2|L_u^*v_s|^2\\
\label{eq:|dt|}&=&-9(\d_t\lam)^2\big|\mu\circ u\big|^2-4\lam^2|L_u^*v_s|^2,
\end{eqnarray} 
\begin{equation}\label{eq:|ds|}-6\d_s(\lam^2)\lan\mu\circ u,L_u^*Jv_s\ran\geq\ldots =-9(\d_s\lam)^2\big|\mu\circ u\big|^2-4\lam^2|L_u^*Jv_s|^2.\end{equation}
Furthermore, a short calculation shows that 
\[-9|d\lam|^2+\La(\lam^2)=-\lam^4\big(|d(\lam^{-1})|^2+\La(\lam^{-2})\big)\geq-C\lam^4.\] 
Combining the estimates (\ref{eq:III}), (\ref{eq:|dt|}) and (\ref{eq:|ds|}), we get 
\begin{equation}\label{eq:La e w lam 4}\La \wt e_w\geq2\lam^4\big(-c^2|\mu\circ u|^2+|L_u\mu\circ u|^2\big)-b\wt e_w^2.\end{equation}
We fix a number $\de>0$ as in Lemma \ref{le:c c}, depending on $c$. Let $z\in\Si$ be a point. If $|\mu\circ u(z)|\leq \de$ then inequality (\ref{eq:La e w lam 4}) implies $\La \wt e_w(z)\geq -b\wt e_w^2(z).$ Suppose now that $|\mu\circ u(z)|>\de$. Then inequality (\ref{eq:La e w lam 4}) implies that 
\[\La \wt e_w(z)\geq -2c^2\de^{-2}\lam^4|\mu\circ u(z)|^4-b\wt e_w^2(z)\geq -C'\wt e_w^2(z),\]
where $C':=2c^2\de^{-2}+b$. This proves Claim \ref{claim:La e w lam}. Lemma \ref{le:e w} follows from Claim \ref{claim:La e w lam} and Lemma \ref{le:mean}. 
\end{proof}
\subsection{Proof of optimal decay}\label{subsec:proof decay}
\begin{proof}[Proof of Theorem \ref{thm:decay}]\setcounter{claim}{0} Assume that the hypothesis (H) of Section \ref{sec:main} is satisfied. Let $\Si,\om_\Si$ and $w$ be as in the hypothesis of Theorem \ref{thm:decay}. We may assume \Wlog that $a=1$, \ie identifying $S^1\iso\R/\Z$, 
\[\Si=\big\{s+it\in \C\,\big|\,s\geq0\big\}/i\Z\iso[0,\infty)\x S^1.\]
A standard argument as in the proof of Proposition D.2 in \cite{Zi} shows that $w$ is gauge equivalent to a smooth vortex, via a locally $W^{2,p}$-gauge transformation defined on $(0,\infty)\x S^1$. (See also \cite{CGMS}.) Hence we may assume \Wlog that $w$ is smooth. We define the function $E:[0,\infty)\to\R$ by $E(s):=E\big(w,[s,\infty)\x S^1\big)$. Let $\eps>0$.
\begin{claim}\label{claim:E s} There exists a number $s_0\geq1$ such that for every $s\geq s_0$ we have
  \begin{equation}
    \label{eq:E s}\frac d{ds}E(s)\leq-(4\pi-\eps)E(s).
  \end{equation}
\end{claim}
\noindent{\bf Proof of Claim \ref{claim:E s}.} By the discussion at the beginning of Subsection \ref{subsec:en act bound} we may assume \Wlog that $P=\Si\x\G$ is the trivial bundle, and we may view $u$ as a map from $\Si$ to $M$. Furthermore, we define $(\Phi,\Psi):\Si\to\g$ by formula (\ref{eq:A z g}). It follows from hypothesis (H) and Lemma \ref{le:c c} that there exists a number $\de_0>0$ such that $K:= \mu^{-1}(\bar B_{\de_0})$ is compact, $\G$ acts freely on $K$, and 
\begin{equation}\label{eq:m mu} m:=m_{\mu^{-1}(\bar B_{\de_0})}\geq\sqrt{\frac{4\pi-\eps}{4\pi}}m_{\mu^{-1}(0)},
\end{equation}
where $m_X$ is as in (\ref{eq:m X}), for $X\sub M$. We fix a number $\de>0$ as in Theorem \ref{thm:isoperi} (Sharp isoperimetric inequality), corresponding to the set $K:=\mu^{-1}(\bar B_{\de_0})$ and the constant $c:=\frac1{4\pi-\eps}$. Shrinking $\de$ we may assume that it also satisfies the condition of Proposition \ref{prop:en act} (Energy action identity) with the same $K$. Using the first inequality in (\ref{eq:lam sup}) with $a=1$, Lemma \ref{le:e w} implies that there exists a number $s_0>0$ such that for $s\geq s_0$ 
\[u(s,t)\in K,\,\forall t\in S^1,\qquad \bar\ell(G u(s,\cdot))\leq \int_0^1|\d_tu+L_u\Psi|(s,t)\,dt<\de.\] 
Hence by the assertion of Theorem \ref{thm:isoperi} with $p:=2$, we have for $s\geq s_0$,
\begin{equation}\label{eq:A 4 pi eps} \big|\A\big((u,\Psi)(s,\cdot)\big)\big|\leq \frac1{4\pi-\eps}\big|\big|\big(\d_tu+L_u\Psi\big)(s,\cdot)\big|\big|_2^2+\frac\pi{m^2}\big|\big|\mu\circ u(s,\cdot)\big|\big|_2^2.
\end{equation}
We fix two numbers $s'\geq s\geq s_0$. By the assertion of Proposition \ref{prop:en act} the loops $u(s,\cdot)$ and $u(s',\cdot)$ are admissible, and 
\begin{equation}
  \label{eq:E w A A}E\big(w,[s,s']\x S^1\big)=-\A\big((u,\Psi)(s',\cdot)\big)+\A\big((u,\Psi)(s,\cdot)\big).
\end{equation}
It follows from inequality (\ref{eq:A 4 pi eps}), the first inequality in (\ref{eq:lam sup}) and Lemma \ref{le:e w} that $\big|\A\big((u,\Psi)(s',\cdot)\big)\big|\to0$, as $s'\to\infty$. Combining this with equality (\ref{eq:E w A A}), inequalities (\ref{eq:A 4 pi eps}) and (\ref{eq:m mu}), and the fact $|\phi|_\op\leq|\phi|$, for $\phi\in\g^*$, we get
\begin{eqnarray}\nn E(s)&=&E\big(w,[s,\infty)\x S^1\big)\\
\nn&=&\A\big((u,\Psi)(s,\cdot)\big)\\
\nn&\leq&\frac1{4\pi-\eps}\int_0^1\left(|\d_tu+L_u\Psi|^2+\frac{4\pi^2}{m_{\mu^{-1}(0)}^2}|\mu\circ u|^2\right)\bigg|_{(s,t)}\,dt\\
\nn&\leq&-\frac1{4\pi-\eps}\frac d{ds}\int_s^\infty\int_0^1\wt e_w(s,t)\,dt\,ds\\
\nn&=&-\frac1{4\pi-\eps}\frac d{ds}E(s).
\end{eqnarray}
Here in the fourth step we used the first inequality in (\ref{eq:lam sup}) with $a=1$ and the definition (\ref{eq:wt e w}) of $\wt e_w$. Claim \ref{claim:E s} follows from this.

By Claim \ref{claim:E s} the derivative of the function $[s_0,\infty)\ni s\mapsto E(s)e^{(4\pi-\eps)s}$ is non-positive, and hence this function is non-increasing. Combining this with Lemma \ref{le:e w}, and recalling the definitions (\ref{eq:e w E}) and (\ref{eq:wt e w}) of $e_w$ and $\wt e_w$, it follows that there exists a constant $C$ such that inequality (\ref{eq:d A u lam}) holds. This completes the proof of Theorem \ref{thm:decay}. \end{proof}
\appendix 
\section{Inequality for the holonomy of a connection}\label{sec:holon}
Let $\G$ be a compact Lie group, $X$ be a (smooth) manifold without boundary, $\pi:P\to X$ be a (smooth) principal $\G$-bundle over $X$, $x\in C^\infty(S^1,X)$ be a loop, and let $p_0\in \pi^{-1}(x(0+\Z))$. Here we identify $S^1\iso\R/\Z$. We denote by $\A(P)$ the space of (smooth) connection one-forms on $P$, and fix $A\in\A(P)$. Recall that the holonomy $h\in \G$ of $A$ around the loop $x$, with base point $p_0$, is defined by the condition $p(1)=p_0h,$ where $p\in C^\infty([0,1],P)$ is the unique horizontal lift of $x$ starting at the point $p_0$. This means that $A_p\dot p=0$ and $p(0)=p_0.$ We choose a Riemannian metric $g_X$ on $X$ and a distance function $d$ on $\G$ that is induced by some Riemannian metric $\lan\cdot,\cdot\ran_\G$ on $\G$. 
\begin{prop}\label{prop:holon} Let $\G,d,X,g_X$ and $P$ be as above, and let $K\sub X$ be a compact subset. Then there exists a constant $C$ satisfying the following condition. If $A\in\A(P)$, $x\in C^\infty(S^1,K)$ is a loop of length $\ell(x)$ less than the injectivity radius $\iota_K$ of $K$ in $X$, and $p_0\in \pi^{-1}(x(0+\Z))$, then the holonomy $h\in\G$ of $A$ around $x$, with base point $p_0$, satisfies the inequality
  \begin{equation}
    \label{eq:d 1 h}d(\one,h)\leq C||F_A||_{L^\infty(X)}\ell(x)^2.
  \end{equation}
\end{prop}
\begin{proof}[Proof of Proposition \ref{prop:holon}]\setcounter{claim}{0} Let $\G,d,P,X,g_X$ and $K$ be as in the hypothesis. Since $\G$ is compact, we may assume \Wlog that $\lan\cdot,\cdot\ran_\G$ is induced by an invariant inner product on $\g$. Let $A\in\A(P)$, $x\in C^\infty(S^1,K)$ be a loop of length less than $\iota_K$, and let $p_0\in \pi^{-1}(x(0+\Z))$. For $t\in \R/\Z\iso S^1$ we define $v(t)\in T_{x_0}X$ to be the unique vector such that 
  \begin{equation}
    \label{eq:x 0 v}\exp_{x_0}v(t)=x(t),\quad ||v(t)||<\iota_K/2.
  \end{equation}
We define $u:[0,1]\x[0,1]\to X$ by $u(s,t):=\exp_{x_0}sv(t)$. There exists a unique smooth map $p:[0,1]\x[0,1]\to P$ satisfying 
\begin{equation}\pi\circ p=u,\qquad p(0,t)\const p_0,\,\forall t\in[0,1],\qquad \label{eq:A ds p}A(\d_sp)\const0.
\end{equation}
We define $\Psi:=A(\d_t p):[0,1]\x[0,1]\to\g$. Let $g:[0,1]\to \G$ be the unique smooth solution of the ordinary differential equation
\begin{equation}
  \label{eq:dot g}\dot g=-\Psi(1,\cdot)g,\qquad g(0)=\one.
\end{equation}
Inequality (\ref{eq:d 1 h}) will now be a consequence of the following two claims. The proof of the first one is straight-forward. 
\begin{claim}\label{claim:holon} $g(1)\in\G$ equals the holonomy of $A$ around $x$, with base point $p_0$.
\end{claim}
\begin{claim}\label{claim:d g 1} There exists a constant $C$ depending only on the Riemannian manifold $(X,g_X)$ and on the compact subset $K\sub X$, such that
\[d(\one,g(1))\leq C\ell(x)^2||F_A||_{L^\infty(X)}.\] 
\end{claim}
\noindent{\bf Proof of Claim \ref{claim:d g 1}.} It follows from (\ref{eq:dot g}) that 
  \begin{equation}\label{eq:d one} d(\one,g(1))\leq \int_0^1|\Psi(1,t)|\,dt\leq\int_0^1\int_0^1|\d_s\Psi(s,t)|\,ds\,dt.\end{equation}
Furthermore, by the definition of $\Psi$ and the last equality in (\ref{eq:A ds p}) we have 
\[p^*F_A=d(p^*A)+(1/2)[p^*A\wedge p^*A]=\d_s\Psi\,ds\wedge dt+0,\]
and hence 
\begin{equation}\label{eq:d s Psi F}|\d_s\Psi|\leq|\d_s u|\,|\d_t u|\,\big(|F_A|\circ u\big).\end{equation} 
Furthermore, by the definition of $u$ and (\ref{eq:x 0 v}) we have 
\begin{equation}\label{eq:d s u v}|\d_su|=|v|\leq\ell(x)/2,\end{equation}
\begin{equation}\label{eq:d t u C 2}|\d_t u(s,t)|=\big|d\exp_{x_0}(sv(t))s\dot v(t)\big|\leq C_1|\dot v(t)|\leq C_2|\dot x(t)|.
\end{equation}
Here $C_1$ and $C_2$ are constants depending only on $X,g_X$ and $K$, and in the third step we used the fact that $\dot v(t)=d\exp_{x_0}(v(t))\big)^{-1}\dot x(t)$. Inserting inequalities (\ref{eq:d s u v}) and (\ref{eq:d t u C 2}) into (\ref{eq:d s Psi F}), we obtain
\[|\d_s\Psi|\leq (C_2/2)\ell(x)|\dot x|\,||F_A||_{L^\infty(X)}.\]
Claim \ref{claim:d g 1} follows by plugging this into (\ref{eq:d one}). This completes the proof of Proposition \ref{prop:holon}.\end{proof}

\begin{thebibliography}{99}
\bibitem[CGMS]{CGMS} K. Cieliebak, A. R. Gaio, I. Mundet i Riera and D. A. Salamon, \emph{The symplectic vortex equations and invariants of Hamiltonian group actions}, J. Symplectic Geom. {\bf 1} (2002), {\bf no. 3}, 543-645.
\bibitem[CGS]{CGS} K. Cieliebak, A. R. Gaio and D. A. Salamon, \emph{$J$-holomorphic curves, moment maps, and invariants of Hamiltonian group actions},  Internat. Math. Res. Notices  2000,  {\bf no. 16}, 831--882.
\bibitem[GS]{GS} A. R. Gaio and D. A. Salamon, \emph{Gromov-Witten invariants of symplectic quotients and adiabatic limits}, J. Symplectic Geom. {\bf 3} (2005), {\bf no. 1}, 55-159.
\bibitem[JT]{JT} A. Jaffe and C. Taubes, \emph{Vortices and monopoles. Structure of static gauge theories}, Progress in Physics, {\bf 2}. Birkh\"auser, Boston, Mass., 1980.
\bibitem[MS]{MS} D. McDuff and D. A. Salamon, \emph{J-Holomorphic Curves and symplectic topology}, AMS Colloquium Publications, {\bf 52}, Providence, RI, 2004.
\bibitem[Mu1]{Mu1} I. Mundet i Riera, \emph{Yang-Mills-Higgs theory for symplectic fibrations}, Ph.D. thesis, Universidad Autonoma de Madrid, April 1999.
\bibitem[Mu2]{Mu2} I. Mundet i Riera, \emph{Hamiltonian Gromov-Witten invariants}, Topology {\bf 42}, {\bf no. 3} (2003), 525-553.
\bibitem[Zi]{Zi} F. Ziltener, \emph{Symplectic Vortices on the Complex Plane and Quantum Cohomology}, Ph.D. thesis, ETH Z\"urich, May 2006.
\end{thebibliography}
\end{document}